\documentclass[oneside]{amsart}
\usepackage{amsfonts}
\usepackage{amssymb}
\usepackage{amsxtra}
\usepackage{amstext}
\usepackage[english]{babel}

\newcommand{\N}{\mathbb{N}}
\newcommand{\E}{\mathbb{E}}
\newcommand{\Z}{\mathbb{Z}}

\newcommand{\pp}{\mathbb{P}}

\newcommand{\kF}{\mathcal{F}}

\newtheorem {lem} {Lemma} [section]

\newtheorem {theo} {Theorem} [section]

\title[A note on multi-type cookie random walk on integers]
      {A note on multi-type cookie random walk on integers}
\author{Bruno Schapira}

\begin{document}

\begin{abstract} We consider a random walk on integers where at the first visits to a site the walker gets a positive drift, but where after a certain number of visits the walker gets a negative drift. We prove that the walker is almost surely transient to the left with positive speed. This is a variant of a model studied by Zerner, Kosygina and Zerner, and Basdevant and Singh. 
\end{abstract}

\maketitle
\let\languagename\relax

\noindent \textbf{Key words:} Excited random walk, Cookie random walk, Self-interacting random walk.

\bigskip

\noindent \textbf{A.M.S. classification:} 60F05; 60F20.

\section{Introduction} 

This note is meant as a small complement to the papers of Zerner \cite{Zer,Zer2}, Kosygina and Zerner \cite{KZer}, and Basdevant and Singh \cite{BaS}. In fact our goal is to show how the tools and results of these papers allow to study another case of the same general model. The problem addressed in \cite{Zer} and \cite{BaS}, which is a generalization of the model introduced by  Benjamini and Wilson \cite{BW}, is roughly the following: let $(p_i)_{i\ge 1}$ be a sequence of real numbers in $[1/2,1)$, such that $p_i=1/2$ for $i$ large enough, let say $i>M$. Then consider the random walk $(X_n)_{n\ge 0}$ on $\Z$ with transition probabilities given by 
$$\pp[X_{n+1}=x+1 \mid X_n =x]=1-\pp[X_{n+1}=x-1 \mid X_n =x]= p_i,$$ if $|\{j\le n \mid X_j =x\}|=i$. 
The question is to describe the asymptotic behavior of this random walk. Zerner \cite{Zer} has proved that it is either recurrent if $\sum_i(2p_i-1)\le 1$, or converges toward $+\infty$, if $\sum_i (2p_i-1)>1$. Then Basdevant and Singh \cite{BaS} have proved that it converges toward $+\infty$ with positive speed if, and only if, $\sum_i(2p_i-1) > 2$. Moreover they give in \cite{BaS2} the exact rate of growth in the null speed regime. In fact Zerner has obtained his results in the more general setting where the sequence $(p_i)_{i\ge 1}$ is random and non necessarily equal to $1/2$ for $i$ large enough. Furthermore he has extended these results in higher dimension in \cite{Zer2}, and recently to the case $p_i\in (0,1)$ with Kosygina \cite{KZer}.

Here we study a variant where $p_i= q \in (0,1/2)$ for all $i>M$. In this case we prove that the random walk always converges toward $-\infty$ with non zero speed. Although this result may seem easier, our proof is not really shorter and we have not found any direct argument to prove it. Apart from its possible own interest, this result may help to give intuition in other models, like 
the generalized reinforced random walks studied in \cite{RS} for instance.

Such problem belongs to the category of excited or reinforced random walks, where transition densities depend on the past trajectory. This subject is treated in the survey of Pemantle \cite{Pem} where the reader can learn in particular about recent developments on vertex or edge reinforced random walks on graphs (mainly $\Z^d$ or trees), and find many other references.

\section{Notation}
We will try to follow as much as possible notations of \cite{Zer} and \cite{BaS}.

An environment is a sequence $w=(w^x(i))_{x\in \Z, i\ge 1}$, where $w^x(i)$ represents the probability for the random walk to go from $x$ to $x+1$ after $i$ visits to $x$. So the initial environment $w_0$ is defined by $w_0^x(i)=p_i$ for all $x\in \Z$ and $i \ge 1$. In particular we take it deterministic. We will sometimes denote for $n\ge 0$ by $w_n$ the (random) environment as seen after $n$ steps. That is $w_n^x(i) = p_j$ where $j$ is equal to $i$ minus the number of visits to $x$ up to time $n-1$. If $w$ and $w'$ are two environments, we note $w\le w'$ if for all $x$ and all $i$, $w^x(i)\le (w')^x(i)$.

We denote by $(\kF_n)_{n\ge 0}$ the natural filtration of $(X_n)_{n\ge 0}$. Thus $w_n$ is $\kF_n$-adapted for all $n\ge 0$. For $x\in \Z$ and $w$ any environment we denote by $\E_{x,w}$ the law of the random walk starting from $x$ with initial environment $w$. We will sometime forget to precise $x$ or $w$ if there is no possible confusion. With these notations our result can be restated as follows: 
\begin{theo}
\label{theo}
For any $(p_i)_{i=1,\dots,M} \in (0,1)^M$ and any $q \in (0,1/2)$, 
$$\pp_{0,w_0}[X_n \to -\infty] = 1,$$
where $w_0^x(i)=p_i$ if $i\le M$ and $w_0^x(i)=q$ otherwise, for any $x\in \Z$. Moreover $\pp_{0,w_0}$ a.s. there exists $v<0$ such that  
$$\lim_{n\to +\infty} \frac{X_n}{n} =v.$$
\end{theo}

\section{Proof}
The first thing is to obtain a zero-one law. For this  
we need a monotonicity result proved by Zerner \cite{Zer}. For $n\in \Z$, we set $T_n=\inf \{t\ge 0 \mid X_t =n\}$.  

\begin{lem}[\cite{Zer} Lemma 15]
\label{monotone}
Let $-\infty\le x\le 0\le z \le +\infty$ and $t\in \N \cup \{\infty\}$. If $w\le w'$, then  
$$\pp_{0,w}[T_z\le T_x \wedge t] \le \pp_{0,w'}[T_z\le T_x\wedge t].$$ 
\end{lem} 
This lemma was stated in \cite{Zer} for environments with only non negative drift, but Zerner's proof applies for general $w$ and $w'$.

In particular, since $\{X_n \to -\infty\}=\lim_{n\to +\infty} \{T_n=+\infty\}$, we get that $w\le w'$ implies $\pp_{w'}[X_n \to -\infty] \le \pp_{w}[X_n \to -\infty]$.

So we need only to prove Theorem \ref{theo} for $w_0$ of the form $w_0^x(i)= p>1/2$ for all $x\in \Z$ and $i\le M$, and $w_0^x(i) =  q < 1/2$ for all $x\in \Z$ and $i>M$. 
For such environment the monotonicity result implies directly a zero-one law: 
\begin{lem} 
For $w_0$ as above, we have 
$$\pp_{w_0}[X_n\to + \infty] = 0 \textrm{ or } 1,$$
and the same statement occurs with $-\infty$ instead of $+\infty$.
\end{lem}
\begin{proof} Let $A=\{X_n \to +\infty\}$. Consider the bounded martingale 
$$M_n= \E_{0,w_0}[1_A \mid \kF_n].$$
Then a general result implies that a.s. $M_n \to 1_A$ as $n\to +\infty$. 
On the other hand if $w_n$ is the environment after $n$ steps, then clearly $w_n^{X_n+x}(i) \le w_0^x(i)$ for any $x\in \Z$ and $i\ge 1$. So by the Markov property of the couple $(X_n,w_n)_{n\ge 0}$ and Lemma \ref{monotone} we have 
$$M_n=\E_{X_n,w_n}[1_A] \le \E_{0,w_0}[1_A]=M_0 \quad \forall n.$$
So $M_n$ is constant, which proves the lemma. 
\end{proof}

The rest of the proof is divided in two parts. First we prove that $\pp[X_n\to + \infty]=0$. Then we prove that the random walk converges toward $-\infty$ with negative speed. 

\subsection{The random walk cannot converge toward $+\infty$}

This will be proved by following a method used in \cite{BaS}: we assume by contradiction that $X_n \to +\infty$, and in particular that $T_n <+ \infty$ for any $n\ge 0$. Next we need to introduce some additional notation. For $n\ge 1$ and $x\in \Z$ let 
$$U_n^x= |\{0\le j < T_n \mid X_j = x, X_{j+1}=x-1\}|.$$  
Let $(B_i)_{i\ge 1}$ be a sequence of independent Bernouilli random variables with parameter $p$ if $i\le M$ and $q$ if $i>M$. For $j\in \N$ define 
$$k_j=\min \{ k \ge 1 \mid |\{1\le i\le k,\ B_i = 1\}|=j+1\}$$
and 
$$A_j = k_j-j-1.$$ 
Then as in \cite{BaS} Lemma 2.1, we have the equality in law: 
$$A_j = A_{M-1} + \xi_1+ \dots + \xi_{j-M+1},$$
where $(\xi_i)_{i\ge 1}$ is a sequence of i.i.d. geometrical random variables with parameter $q$.
We define now the Markov chain $(Z_n)_{n\ge 0}$ by 
$$\pp[Z_{n+1}= k \mid Z_n = j] = \pp[A_j = k].$$
This is in fact a branching process with migration. Let us recall some facts: first for each $n$, $(U^n_n,\dots,U^0_n)$ has the same law as $(Z_0,\dots,Z_n)$ (see \cite{BaS} Proposition 2.2). In particular since $U^0_\infty$ is a.s. finite, this implies that $(Z_n)_{n\ge 0}$ is irreducible positive recurrent, and converges in law toward some random variable $Z_\infty$, which is independent of the starting point of $(Z_n)_{n\ge 0}$. 
The random generating function of $Z_\infty$ is defined by: 
$$G(s) = \E[s^{Z_\infty}]=\sum_{k\ge 0} \pp[Z_\infty = k] s^k \quad \forall s \le 1.$$
Remember that the law of $Z_\infty$ is the stationary measure of $(Z_n)_{n\ge 0}$. Thus $G(s)=\E[\E_{Z_\infty}[s^{Z_1}]]$, which as in \cite{BaS} Lemma 3.4, leads to the equation 
\begin{eqnarray}
\label{equationG}
1-G\left(\frac{q}{1-(1-q)s}\right)=a(s)(1-G(s))  + b(s),
\end{eqnarray}
where 
$$a(s) =\frac{1}{\E[s^{A_{M-1}}] \left(\frac{1-(1-q)s}{q}\right)^{M-1}}$$
and 
\begin{eqnarray*}
b(s)=1-a(s)+a(s)\sum_{k=0}^{M-2} \pp[Z_\infty = k]\left\{\E[s^{A_k}]-\E[s^{A_{M-1}}]\left(\frac{1-(1-q)s}{q}\right)^{M-1-k}\right\}.
\end{eqnarray*} 
Now observe that all the derivatives of $a$ and $b$ in $1$ are finite. So \eqref{equationG} implies that all the derivatives of $G$ in $1$ are also finite, i.e. all the moments of $Z_\infty$ are finite. Since $Z_\infty$ is stochastically decreasing with $q$, this would imply the same result when $q=1/2$. But this is in contradiction with a general result of Kaverin \cite{Kav} (see also \cite{BaS} Lemma 3.3). Thus $\pp[X_n \to + \infty]=0$ for all $q<1/2$ and all $p<1$.

\subsection{The random walk is transient to the left with negative speed}
For $x\in \Z$ and $n\in \N\cup\{+\infty\}$, define the drift accumulated in $x$ at time $n$ by $$D^x_n=\sum_{i=1}^{|\{j\le n \mid X_j=x\}|}(2p_i-1),$$ and 
the total drift accumulated in the negative part by $$D^-_n=\sum_{x\le 0} D^x_n.$$ 
\begin{lem} \label{dinfini0} We have
$$\E[D^0_\infty]=-1.$$
\end{lem}
\begin{proof} The proof is adapted from \cite{Zer} Lemma 6 and Lemma 11, but is slightly different. Thus let us give some details. 
The basic tool for proving the lemma  is the martingale defined by 
$$R_n := X_n^- - D^-_n \quad \forall n,$$ 
where $X_n^-=\sum_{m\le n-1} (X_{m+1}-X_m)1_{\{X_m \le 0\}}$.  
Observe that if 
$$V_n= \sum_{i\le n-1} 1_{\{X_i=0,X_{i+1}=1\}} \quad \forall n\ge 1,$$ 
then 
$$X_n^- = X_n \wedge 0 + V_n.$$
The fact that $(R_n)_{n\ge 0}$ is a martingale implies that for any $t\ge 0$
$$\E[X_{T_{-K}\wedge t}]=\E[D^-_{T_{-K}\wedge t}]-\E[V_{T_{-K}\wedge t}].$$
Now remark that since $X_n$ a.s. does not tend to $+\infty$, $T_{-K}$ is a.s. finite for any $K\ge 0$.   
So 
$$\lim_{t\to +\infty} \E[X_{T_{-K}\wedge t}\wedge 0] =-K.$$
Moreover the monotone convergence theorem implies that 
$$\lim_{t\to +\infty} \E[V_{T_{-K}\wedge t}] = \E[V_{T_{-K}}].$$
In the same way, observe that the positive and negative parts of $D^-_n$ are increasing with $n$, and that its positive part is bounded for $n\le T_{-K}$, if $K$ is fixed. So again the monotone convergence theorem gives
$$\lim_{t\to +\infty}\E[D^-_{T_{-K}\wedge t}]=\E[D^-_{T_{-K}}],$$
which leads to the equation
\begin{eqnarray}
\label{equationdrift} 
\E[D^-_{T_{-K}}]-\E[V_{T_{-K}}]=-K.
\end{eqnarray}
Now assume that $\E[V_\infty]=+\infty$. Then, since up to multiplicative and additive constants $V_n$ is lower than $D^0_n$, we would have 
$\E[D_\infty^0]=-\infty$. We claim that this would contradict \eqref{equationdrift}. 
Indeed for $x\le 0$, define recursively $\tau^x(m)$ for $m\ge 0$ by $\tau^x(0)=0$ and 
$$\tau^x(m+1)=\inf\{t > \tau^x(m) \mid X_t \le x\}.$$ 
Then observe that the law of $(X_{\tau^x(m)})_{m\ge 0}$ is independent of $(w_n^y)_{n\ge 0, y>x}$. This shows that for any $k\in \N \cup\{+\infty\}$, the law of $D_{T_{x-k}}^x$ is independent of $x\le 0$ under $\E_0$. Thus 
\begin{eqnarray}
\label{equationdrift2}
\E_{0,w_0}[D^-_{T_{-K}}]=\sum_{x=-K}^0 \E_{x,w_{T_x}}[D^x_{T_{-K-x}}]=\sum_{x=-K}^0 \E_{0,w_0}[D^0_{T_{-K-x}}].
\end{eqnarray}
But again a monotone convergence argument yields to 
$$\lim_{K\to +\infty}\E[D^0_{T_{-K}}] = \E[D_\infty^0].$$
Thus we would have 
$$\lim_{K\to +\infty} \frac{\E[D^-_{T_{-K}}]}{K} = -\infty,$$
contradicting \eqref{equationdrift}. So $\E[V_\infty]<+\infty$. Then the lemma follows from \eqref{equationdrift} and \eqref{equationdrift2}.  
\end{proof}

We can now conclude the proof of our theorem. First the random walk cannot be recurrent. Otherwise a.s. $T_n$ would be finite for any $n\in \Z$. Thus a.s. $X_n$ would be equal to $0$ infinitely often, contradicting the above lemma. So a.s. $X_n \to -\infty$. Moreover as in \cite{BaS}, one can see that the random walk has a negative speed if and only if $\E[V_\infty]<+\infty$. Thus again the above lemma implies that it has well a negative speed.

\vspace{0.2cm}
\noindent  D\'epartement de Math\'ematiques, B\^at. 425, Universit\'e Paris-Sud, F-91405 Orsay
cedex, France.

e-mail: bruno.schapira@math.u-psud.fr

\end{document}